\documentclass[10pt,reqno]{amsart} 
\usepackage{preamble}
\title{Labelling Schubert intersections in the Grassmanian}
\author{Noah White}
\address{Noah White, Mathematics Department, UCLA, Los Angeles, CA 90095, USA}
\email{noah@math.ucla.edu}

\begin{document}
\begin{abstract}
  Points in the intersection of Schubert varieties are counted by various combinatorial objects, for example standard tableaux. This paper consider the problem of producing a natural labelling of intersection points by these combinatorial objects. When the Schubert varieties are being taken with respect to flags osculating at real points, several different methods have appeared implicitly in the literature (specifically in work of Mukhin-Tarasov-Varchenko, Speyer and Marcus). In this paper we show these various methods produce the same labelling and we describe it in an elementary way.
\end{abstract}

\maketitle
\tableofcontents

\section{Introduction}
\label{sec:introduction}

Intersections of Schubert varieties in the Grassmanian provide the prototypical examples of many phenomena in geometry, representation theory and combinatorics. Let \( V \) be a finite dimensional complex vector space and \( \Gr(r,V) \), the \emph{Grassmanian} variety of \( r \)-dimensional subspaces of \( V \). To every partition \( \lambda \) with at most \( r \) rows and \( d = \dim V \) columns, and a full flag \( \Ff: \{0 \} = \Ff_0 \subset \Ff_1 \subset \cdots \subset \Ff_d = V \), we can associate the \emph{Schubert variety} \( \Shvar(\lambda,\Ff) \) of subspaces that meet \( \Ff \) in a way prescribed by \( \lambda \) (see Section~\ref{sec:preliminaries}).

The intersection theory of Schubert varieties is governed by the \emph{Littlewood-Richardson coefficients} \( c^\nu_{\lambda\mu} \). For a generic choices of flags \( \Ff,\Gg,\Hh \) (and appropriate partitions) the intersection
\[ \Shvar(\lambda,\Ff) \cap \Shvar(\mu,\Gg) \cap \Shvar(\nu^\comp,\Hh) \]
is a set of \( c^\nu_{\lambda\mu} \) points. Here \( \nu^\comp \) is the complement partition (see Section~\ref{sec:standard-tableaux}).

The Littlewood-Richardson coefficients have a multitude of different incarnations. They count \emph{Littlewood-Rchardson tableaux} and \emph{dual equivalence classes}. They count the multiplicity of the irreducible \( \gl_r \)-module \( L(\nu) \) in the tensor product \( L(\lambda) \otimes L(\mu) \), and they are the structure constants of the ring of symmetric functions for the basis of Schur polynomials.

More generally one can consider intersections of the form
\begin{equation}
  \label{eq:schubert-intersection}
  \Shvar(\lambda^{(1)},\Ff^{(1)}) \cap \Shvar(\lambda^{(2)},\Ff^{(2)}) \cap \cdots \cap \Shvar(\lambda^{(n)},
  \Ff^{(n)}) \cap \Shvar(\mu^\comp,\Ff^{(\infty)})  
\end{equation}
for a sequence of partitions \( \blambda = (\lambda^{(1)},\lambda^{(2)},\ldots,\lambda^{(n)}) \). For generic choices of flags (and appropriate choices of partitions) this intersection will again be finite and counted by Littlewood-Richardson coefficients \( c^\mu_{\blambda} \). This paper will concentrate on the special case when \( \lambda^{(i)} = (1) = \square \) for all \( i \). In that case, when \( \mu \) is a partition of \( n \), \( c^\mu_{\square,\square,\ldots,\square} \) is equal to the number of standard tableaux of shape \( \mu \) (which we will often think of as chains of shapes with terminal shape being \( \mu \)).

\subsection{Labelling Schubert intersections}
\label{sec:labell-schub-inters}

The obvious question which arises is whether there is a canonical bijection that realises this coincidence of numbers. In this form, the answer to the question is no. As one varies the flags, there can be monodromy, often the entire symmetric group will act on the fibre. However, implicit in work over the past decade are a number of ways to realise bijections between the above intersections and the set of standard tableaux when the flags are chosen to be osculating flags at real numbers. In short, we aim to show these bijections are actually the same and describe this bijection in elementary geometric terms.

We start with a rational normal curve \( \PP^1 \longrightarrow \PP V \). For any \( z \in \PP^1 \), let \( \Ff(z) \) be the full flag of \emph{osculating} subspaces to \( z \in \PP^1 \), that is, \( \Ff(z)_i \) is the unique \( i \)-dimensional subspace of \( \PP V \) having maximal intersection with the curve at \( z \).

The first bijection comes from work of Mukhin, Tarasov and Varchenko~\cite{Mukhin:2009et}. The authors show that the ring of functions on the scheme theoretic intersection \( \Shvar(\blambda,\mu^\comp;z,\infty) \) given in~(\ref{eq:schubert-intersection}), with \( \Ff^{(i)} = \Ff(z_i) \) for a tuple of distinct complex numbers \( z = (z_1,z_2,\ldots,z_n) \), is isomorphic to a certain commutative algebra of operators on
\[ L(\blambda)^{\sing}_\mu = \left[ L(\lambda^{(1)}) \otimes L(\lambda^{(2)}) \otimes \cdots \otimes L(\lambda^{(n)}) \right]^{\sing}_\mu, \]
the space of singular (i.e. highest weight) vectors of weight \( \mu \) in a tensor product of irreducible \( \gl_r \)-representations. If \( z \) is taken to be real such that \( z_1<z_2<\cdots <z_n \), then one can take a limit as \( z \to \infty \) (in a specified way), this commutative algebra tends to the algebra of \emph{Jucys-Murphy operators}. When \( \blambda = \square^n \) the spectrum of the Jucys-Murphy operators is in canonical correspondence with the set of standard tableaux of shape \( \mu \).

The second bijection comes from the work of Speyer~\cite{Speyer:2014gg}. The author constructs an explicit bijection between the points in~(\ref{eq:schubert-intersection}) and certain \emph{cylindrical growth diagrams}. Again, there is a natural way to place these objects in bijection with standard tableaux of shape \( \mu \).

\begin{Theorem}
  \label{thm:labelling-theorem}
  The bijections between \( \Shvar(\square^n,\mu^\comp;z,\infty) \) and \( \SYT(\mu) \), the set of standard tableaux of shape \( \mu \), defined by Speyer and Mukhin-Tarasov-Varchenko agree.
\end{Theorem}

\begin{Remark}
  \label{rem:general-intersections}
In this paper we only consider intersections of Schubert varieties for sequences of partitions \( \blambda \) where \( \lambda_i = \square \) for all \( i \), however the question can be asked for any sequence of partitions. The methods of this paper can be applied in this more general setting, one only needs to understand what combinatorial objects should label the intersection points. If we think of a standard tableau as a chain of partitions, one box added at a time, then the correct generalisation would be a chain of partitions with \( \left| \lambda_i \right| \) boxes added at the \( i^{\text{th}} \) step, as well as the data of a dual equivalence class labelling the \( i^{\th} \) inclusion, in such a way that this dual equivalence class is slide equivalent to \( \lambda_i \). 
\end{Remark}

\subsection{An elementary description} 
\label{sec:an-elem-descr}

We give an elementary description of the bijection described in Theorem~\ref{thm:labelling-theorem}. We fix an \( n \)-tuple of real numbers \( z = (z_1,z_2,\ldots,z_n) \) such that \( z_1 < z_2 < \cdots < z_n \). Choose a subspace in the intersection
\[ X \in \Shvar(\square^n,\mu^\comp; z,\infty). \]
This point depends, in particular on \( z_n \in \RR \). We analyse what happens when we take the limit \( X_\infty = \lim_{z_n \to \infty} X \). The limit point exists since the Grassmanian is projective.

\begin{Proposition}
  \label{prp:limit-exists-intro}
The limit point \( X_\infty \in \Gr(r,V) \) is contained in the Schubert cell \( \Shvar^{\circ}(\lambda^\comp;\infty) \) for some partition \( \lambda \) obtained from \( \mu \) by removing a single box.
\end{Proposition}

In particular this shows that \( X_\infty \in \Shvar(\square^{n-1},\lambda^\comp;z_1,\ldots,z_{n-1},\infty) \) and by induction we can associate to \( X_\infty \) a standard \( \lambda \)-tableaux \( S \). We then associate to \( X \) the unique standard \( \mu \)-tableau which is equal to \( S \) after removing the box containing \( n \).

\begin{Theorem}
  \label{thm:alternative-bij}
The above process describes a bijection \( \Shvar(\square^n,\mu^\comp;z,\infty) \longrightarrow \SYT(\mu) \). Furthermore it coincides with the bijection from Theorem~\ref{thm:labelling-theorem}.
\end{Theorem}

\subsection{Method of proof}
\label{sec:method-proof}

We will first show that Speyer's bijection agrees with the bijection described in Section~\ref{sec:an-elem-descr}. To make the connection with Bethe vectors we use an intervening step, the critical points of the master function.

The master function is a certain function whose critical points give rise to Bethe vectors and from which one can also determine the spectrum of the Gaudin Hamiltonians. Mukhin-Tarasov-Varchenko give a map between points in Schubert intersections and critical points. We investigate the compatibility of the process outlined in Section~\ref{sec:an-elem-descr} with this map.

The asymptotics of critical points were investigated by Reshetikhin and Varchenko in~\cite{Reshetikhin:1995vs} and this result was exploited by Marcus~\cite{Marcus:2010vn} to provide a way of labelling critical points by standard tableaux. We will recall this as well as the proofs for the sake of convenience. This description will allow us to more easily compare the MTV and Speyer labellings.

\subsection{Structure of paper}
\label{sec:structure-paper}

In Section~\ref{sec:schubert-intersections} we review the basic definitions of Schubert varieties and outline an elementary procedure which associates a tableau to a point in a Schubert intersection. Sections~\ref{sec:mtv-labelling} and~\ref{sec:speyers-labelling} then define the MTV and Speyer labellings respectively and in Section~\ref{sec:agre-with-elem} a proof is given that the Speyer labelling coincides with the elementary labelling. Finally, in Sections~\ref{sec:algebr-bethe-ansatz} and~\ref{sec:crit-points-schub} the technology of the algebraic Bethe ansatz is introduced as well as Marcus' analysis of it, this is then used to prove that the Speyer and MTV labellings coincide.

\subsection{Acknowledgements}
\label{sec:acknowledgements}

The author thanks Iain Gordon and Arun Ram for helpful conversations and inspiration.

\section{Schubert intersections}
\label{sec:schubert-intersections}

In this section we give some background on intersections of Schubert varieties, the combinatorics of standard tableaux and use this to define an elementary procedure for labelling the points in certain intersections by standard tableaux.

\subsection{Preliminaries}
\label{sec:preliminaries}
Let \( V \) be a \( d \) dimensional vector space and \( \Gr(r,d) \) the \emph{Grassmanian} of \( r \)-planes in \( V \). Fix a rational normal curve \( \PP^1 \longrightarrow \PP V \)  and for \( z \in \PP^1 \) let \( \Ff(z) \) be the full flag of subspaces of \( V \) osculating at \( z \). For concreteness and without loss of generality, we may choose a basis of \( V \) and identify \( V = \CC_{d-1}[x,y] \), the space of homogeneous degree \( d-1 \) polynomials. \( \Ff([a:b]) \) is then the flag
\[ (bx-ay)^{d-1}\CC_0[x,y] \subset (bx-ay)^{d-2}\CC_1[x,y] \subset \cdots \subset \CC_{d-1}[x,y].  \]
We will often identify \( \CC_d[x,y] \) with the space \( \CC_{<d}[u] \) of polynomials of degree less than \( d \), using the map \( x^iy^{d-i} \mapsto u^i \).

Let \( \lambda \) be a partition with at most \( r \) rows and \( d-r \) columns. Given a flag \( \Ff \) of subspaces of \( V \), the \emph{Schubert cell} \( \Shvar^\circ(\lambda;\Ff) \subset \Gr(r,d) \) is the subvariety of subspaces \( X \) such that
\[ \dim(X \cap F_k) = \# \setc{1 \le s \le r}{d-r+1-\lambda_s \le k}. \]
The closure of \( \Shvar^\circ(\lambda;\Ff) \) is denoted \( \Shvar(\lambda;\Ff) \) and is called the \emph{Schubert variety}. It is the subvariety of subspaces \( X \) such that
\[ \dim(X \cap F_k) \ge \# \setc{1 \le s \le r}{d-r+1-\lambda_s \le k}. \]
If \( \Ff = \Ff(z) \) we denote \( \Shvar(\lambda;\Ff(z)) \) by \( \Shvar(\lambda;z) \).

\subsection{Intersections}
\label{sec:intersections}

We will generally be interested in the intersection of \( k \) Schubert varieties. If \( \blambda = (\lambda^{(1)},\lambda^{(2)},\ldots,\lambda^{(k)}) \) is a sequence of partitions and \( z = (z_1,z_2,\ldots,z_k) \) then we will use the notation
\[ \Shvar(\blambda;z) = \bigcap_{i=1}^k \Shvar(\lambda^{(i)};z_i). \]
It is a theorem of Eisenbud and Harris~\cite{Eisenbud:1983fj} that when the \( z_i \) are distinct, the intersection \( \Shvar(\blambda,z) \) has maximum possible codimension \( \left| \blambda \right| = \sum_{i=1}^k \left|\lambda^{(i)}\right| \), where \( \left| \mu \right| \) is the size of the partition \( \mu \). This result has a strengthening as conjectured by Shapiro-Shapiro and proved by Mukhin-Tarasov-Varchenko.

\begin{Theorem}[\cite{Mukhin:2009cf}]
  \label{thm:shaprio-conj}
When \( z = (z_1,z_2,\ldots,z_k) \) is a set of distinct real points and \( \left| \blambda \right| = r(d-r) \), the intersection \( \Shvar(\blambda;z) \) is a reduced union of real points of \( \Gr(r,d) \).
\end{Theorem}

We let \( X_n = \setc{ z \in \CC^n}{z_i \neq z_j} \) be the set of \( n \)-tuples of distinct complex numbers, \( X_n(\RR) \) the set of real points and \( X_{n}^< \subset X_n(\RR) \) the set of \( n \)-tuples of real numbers \( z=(z_1,z_2,\ldots,z_n) \) such that \( z_1<z_2<\ldots <z_n \).

\subsection{The Wronskian}
\label{sec:wronskian}

Let \( g_1(u),g_2(u),\ldots,g_k(u) \) be a collection of polynomials in the variable \( u \). Recall the \emph{Wronskian} is the determinant
\begin{equation*}
  \Wr(g_1,g_2,\ldots,g_k) = \det
  \begin{pmatrix}
    g_1(u)  & g_2(u)  & \cdots & g_k(u)  \\
    g_1'(u) & g_2'(u) & \cdots & g_k'(u) \\
    \vdots  & \vdots  & \ddots & \vdots \\
    g_1^{(n-1)}(u) & g_2^{(n-1)}(u) & \cdots & g_k^{(n-1)}(u) 
  \end{pmatrix}.
\end{equation*}

\begin{Lemma}
  \label{lem:linear-class-wronskian}
Up to a scalar factor, the Wronskian \( \Wr(g_1(u),\ldots,g_k(u)) \), depends only on the subspace of \( \CC[u] \) spanned by the polynomials \(  g_i(u) \), and is zero if and only if the polynomials are linearly dependant.
\end{Lemma}

\begin{proof}
The result follows from the fact that the derivative is a linear operation, and that the determinant can at most be multiplied by a scalar after applying column operations.
\end{proof}

This allows us to define the \emph{Wronskian} \( \Wr\map{\Gr(r,d)}{\PP(\CC[u])} \) by \( \Wr(X) = \Wr(f_1,f_2,\ldots,f_r) \) for any basis \( \{f_1,f_2,\ldots,f_r\} \) of \( X \).

\subsection{Standard tableaux}
\label{sec:standard-tableaux}
\begin{figure}
  \centering
  \begin{tikzpicture}
    \draw[thick] (0,0) rectangle (4,-2.5);
    \draw[step = 0.25,help lines] (0,0) grid (4,-2.5);
    \draw[thick] (2.5,0) -- (2.5,-0.75) -- (2,-0.75) -- (2,-1) -- (1.75,-1) -- (1.75,-1.75) -- (1,-1.75) -- (1,-2) -- (0,-2);
    \node at (1,-1) {\( \mu \)};
    \node at (3,-1.5) {\( \mu^\comp \)};
  \end{tikzpicture}
  \caption{The complementary partition}
  \label{fig:compl-part}
\end{figure}
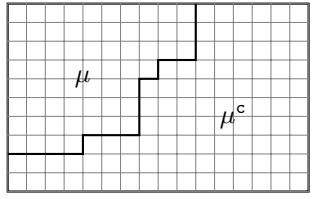
Given a partition \( \mu \) of \( n \), denote by \( \SYT(\mu) \) the set of standard \( \mu \)-tableaux, that is, the set of fillings of the diagram of shape \( \mu \) using \( 1,2,\ldots,n \) with increasing rows and columns. If \( \mu \) has at most \( r \) parts, and \( d-r \) columns, then we denote by \( \mu^\comp \) the partition obtained by embedding \( \mu \) into the top left of a \( r \times (d-r) \) rectangle, taking the complement and rotating the resulting shape by \( 180 \) degrees. See Figure~\ref{fig:compl-part}. I.e.
\[ \mu^\comp  = (d-r-\mu_r,d-r-\mu_{r-1},\ldots,d-r-\mu_1) \]
Repeated application of the Pieri rule and Theorem~\ref{thm:shaprio-conj} imply that when \( z\in X_n(\RR) \) then
\[ \# \Shvar(\square^n,\mu^\comp; z,\infty) = \# \SYT(\mu). \]
For notational simplicity we will set
\[ \Shvar(z)_\mu = \Shvar(\square^n,\mu^\comp; z,\infty) = \Shvar(\square;z_1)\cap \Shvar(\square;z_2) \cap\cdots\cap\Shvar(\square;z_n)\cap \Shvar(\mu^\comp;\infty). \]

\subsection{An explicit bijection}
\label{sec:an-expl-biject}

In this section we will describe an explicit and elementary way to realise the numerical coincidence above as a bijection \( \Shvar(z)_\mu \longrightarrow \SYT(\mu) \). Fix an \( n \)-tuple of real numbers \( z \in X_n^< \) and  fix a subspace \( X \in \Shvar(z)_\mu \). We will associate to it a standard \( \mu \)-tableaux. We will do this by induction on \( n \).

For \( n = 1 \), there is a single point in \( \Shvar(\square,\square^\comp;z,\infty) \) and a single \( \square \)-tableaux so the labelling is completely determined. For \( n > 1 \) let \( X \in \Shvar(z)_\mu \). We fix the first \( n-1 \) marked points at \( z_1,\ldots,z_{n-1} \), the last marked point at \( \infty \) and let the \( n^{\text{th}} \) marked point vary in the following sense.

Let \( z(s) = (z_1,z_2,\ldots,z_{n-1},s) \), then we have the family \( \Omega(z(s))_\mu \) over \( P_n = \RR - \{ z_1, z_2, \ldots, z_{n-1} \} \). Let \( X(s) \) be a local section of \( \Omega(z(s))_\mu \) such that \( X(z_n) = X \). Let
\begin{equation*}
  X_\infty = \lim_{s \rightarrow \infty} p X(s) \in \Gr(r,d),
\end{equation*}
where \( p \) is the projection to \( \Gr(r,d) \).

\begin{Lemma}
  \label{lem:limit-one-partition-smaller}
The limit point \( X_\infty \) exists and is contained in \( \Shvar^\circ(\lambda^\comp,\infty) \) for some partition \( \lambda \subset \mu \) such that \( \left| \lambda \right| + 1 = \left| \mu \right| \). In particular
\begin{equation*}
  X_\infty \in \Shvar(z_1,z_2,\ldots,z_{n-1})_{\lambda}.
\end{equation*}
\end{Lemma}

\begin{proof}
Since the Grassmanian is projective, \( X_\infty \) must exist and since a Schubert variety is a union of smaller Schubert cells,
\begin{equation*}
\Shvar(\mu^\comp;\infty) = \bigsqcup_{\nu \supseteq \mu^\comp} \Shvar^\circ(\nu;\infty),
\end{equation*}
we must have that \( X_\infty \in \Shvar^\circ(\lambda^\comp;\infty) \) for some partition \( \lambda \) such that \( \lambda^\comp \supseteq \mu^\comp \), or equivalently, such that \( \lambda \subseteq \mu \).

Since each of the varieties \( \Shvar(\square;z_i) \) is closed, \( X_\infty \) must lie in the intersection \( \Shvar(z_1,z_2,\ldots,z_{n-1})_\lambda = \Shvar(\square^{n-1},\lambda^\comp;z_1,\ldots,z_{n-1},\infty) \). However this is empty unless \( \left| \lambda \right| \ge n-1 \). Hence either \( \lambda = \mu \) or \( \lambda \) is obtained by removing a single box from \( \mu \). To decide whether \( \abs{\lambda} = n \) or \( n-1 \) we use the Wronskian. By \cite[Lemma~4.2]{Mukhin:2009et} 
\begin{equation*}
\Wr(X(s))(u) = (u-s)\prod_{a=1}^{n-1} (u - z_a).
\end{equation*}
By continuity, \( \Wr(X_\infty) = \prod_{a=1}^{n-1}(u-z_a) \). It is straightforward to see from the definition that if \( Y \in \Shvar^\circ(\nu^\comp;\infty) \) then \( \deg \Wr(Y) = \left| \nu \right| \). We have shown that \( \Wr(X_\infty) = n-1 \) so therefore \( \left| \lambda \right| = n-1 \).
\end{proof}

Summarising, \( X_\infty \in \Shvar(z_1,\ldots,z_{n-1})_\lambda \) and \( \left| \lambda \right| = n-1 \). Thus by induction we can assign a standard \( \lambda \)-tableau to \( X_\infty \). Let this tableau be \( T' \). Let \( \EL(X) \) be the unique, standard \( \mu \)-tableaux which is \( T' \) upon restriction to \( \lambda \subset \mu \).

\begin{Theorem}
  \label{thm:labelling}
The map \( \EL\map{\Shvar(z)_\mu}{\SYT(\mu)}; X \mapsto \EL(X) \) is a bijection.
\end{Theorem}

We will prove this theorem by showing that it coincides with two maps, each of which is known to be a bijection.

\section{The MTV labelling}
\label{sec:mtv-labelling}

Now we describe a labelling of the points of \( \Shvar(z)_\mu \) by standard tableaux implicit in the work of Mukhin-Tarasov-Varchenko using the spectrum of \emph{Bethe algebras} which the aforementioned authors have identified with intersections of Schubert varieties. The idea is that these algebras degenerate to the algebra of Jucys-Murphy operators in a certain limit. 

\subsection{Bethe algebras}
\label{sec:bethe-algebras}

Let \( \glr \) be the general linear Lie algebra and \( \glr[t] \) the \emph{current algebra} of \( \glr \)-valued polynomials in \( t \). For an element \( g \in \glr \), let \( g(u) = \sum_{s \ge 0} gt^su^{-s-1} \) be a formal power series with values in \( U(\glr[t]) \). The differential operator \( \Dd = \det \left( \delta_{ij}\partial_u  - e_{ji}(u) \right) \) is defined by expansion along the first column and can be expended in powers of \( \partial_u \)
\[ \Dd = \sum_{i=0}^r \sum_{s \ge 0} B_{is}u^{-s}\partial_u^{r-i}, \text{ for } B_{is} \in U(\glr[t]). \]
the \emph{universal Bethe algebra} is the subalgebra \( \uBethe \subseteq U(\glr[t]) \)  generated by the coefficients \( B_{is} \). The algebra \( \uBethe \) is a commutative subalgebra of \( U(\glr[t])^{\glr} \), i.e. it commutes with \( U(\glr) \subset U(\glr[t]) \) (see \cite[Proposition~8.2 and~8.3]{Mukhin:2006vt}).

Given a \( \glr \)-module \( M \) and \( z \in \CC \), let \( M(z) \) be the \( \glr[t] \)-module where \( t \) acts as multiplication by \( z \). Let \( L(\lambda) \) be the finite dimensional, irreducible, highest-weight module for \( \glr \), corresponding to the partition \( \lambda \). Then for a tuple of distinct complex numbers \( z \in X_n \), and a tuple of partitions \( \blambda = (\lambda^{(1)},\lambda^{(2)},\ldots,\lambda^{(n)}) \), the algebra \( \uBethe \) acts on
\[ L(\blambda;z)^\sing_\mu = \left[ L(\lambda^{(1)})(z_1) \otimes L(\lambda^{(2)}) \otimes \cdots \otimes L(\lambda^{(n)})(z_n) \right]^{\sing}_{\mu}, \]
the space of highest-weight vectors of weight \( \mu \) in the tensor product. We will be mostly interested in the case \( \blambda = \square^n \) in which case we will denote the above space by \( L(z)^\sing_\mu \). The image of \( \uBethe \) in \( \End(L(z)^\sing_\mu) \) will be denoted \( \uBethe(z)_\mu \).

\subsection{The spectrum and the Gaudin Hamiltonians}
\label{sec:spectr-gaud-hamilt}

The spectrum of the algebra \( \uBethe(z)_\mu \) will be denoted \( \bspec(z)_\mu \). Let \( V = L(\square) \), then \( L(z)_\mu \cong [V^{\otimes n}]^\sing_\mu \) as a \( U(\glr)^{\otimes n} \)-module. Define operators
\[ H_a(z) = \sum_{b \neq a} \frac{(a,b)}{z_a-z_b} \text{ for } 1 \le a \le n, \]
where \( (a,b) \) is the transposition swapping the \( a^\th \) and \( b^{\th} \) tensor factors of \( V^{\otimes n} \). The symmetric group \( \Sn \) acts on \( V^{\otimes n} \) and preserves the subspace \( S_\mu = [V^{\otimes n}]^\sing_\mu  \) which is a copy of the irreducible \( \Sn \)-module corresponding to the partition \( \mu \).

\begin{Theorem}[{\cite[Theorem~3.2 and~Corollary~3.3]{Mukhin:2010ky}}]
  \label{thm:propterties-bethe}
  For generic \( z \in X_n \) (including all real \( z \in X_n(\RR) \)),
  \begin{enumerate}
  \item the algebra \( \uBethe(z)_\mu \) is generated by the operators \( H_a(z) \), and
  \item the algebra \( \uBethe(z)_\mu \) has simple spectrum, that is, \( \bspec(z)_\mu \) is a reduced set of \( \dim S_\mu = \# \SYT(\mu) \) points.
  \end{enumerate}
\end{Theorem}

The coincidence of \( \# \bspec(z)_\mu = \#SYT \) in the above theorem can again be realised using a limiting process. Let \( z \in X_n^<  \). Choose a path \( z(t) = (z_1(t), z_2(t), \ldots,z_n(t)) \in X_n^< \) such that
\begin{enumerate}
\item \( z(1) = z \),
\item \( \lim_{t \to \infty} z_i(t) = \infty \), and
\item \( \lim_{t \to \infty} z_i(t)/z_{i+1}(t) = 0 \).
\end{enumerate}
Then, in this limit \( \lim_{t\to\infty} z_a(t)H_a((z(t))) = L_a = \sum_{b < a} (a,b) \), the \( a^\th \) \emph{Jucys-Murphy operator} on \( S_\mu \). These operators are well known to have simple spectrum on \( S_\mu \), and their spectrum can be canonically identified with \( \SYT(\mu) \) in the following way. If \( v \in S_\mu \) is a joint eigenvector such that \( L_av = c_av \), then we associate to the eigenspace \( \CC v \) the unique standard tableau \( T \) where the box containing \( a \) has content \( c_a \). Thus we obtain a bijection by parallel transport
\[ p:\bspec(z)_\mu \longrightarrow \JMspec_\mu \cong \SYT(\mu) \]
where \( \JMspec_\mu \) is the joint spectrum of the Jucys-Murphy operators on \( S_\mu \). Since the parameter space is contractible, this does not depend on the choice of path.

\subsection{Isomorphism to Schubert intersections}
\label{sec:isom-schub-inters}

Let \( \chi \in \bspec(z)_\mu \) be a closed point and identify this with a functional \( \chi\map{\uBethe(z)_\mu}{\CC} \). Consider the differential operator on \( \CC[u] \) given by
\[ \chi(\Dd) = \sum_{i=0}^r \sum_{s \ge 0} \chi(B_{is})u^{-s}\partial_u^{r-i} \]
and denote the kernel of this operator by \( X_\chi = \ker \chi(\Dd) \). By~\cite[Lemma~5.6]{Mukhin:2004ks} \( X_\chi \) is an \( r \)-dimensional subspace of \( \CC_d[u] \) and moreover \( X_\chi \in \Shvar(z)_\mu \subset \Gr(r,d) \).

\begin{Theorem}[{\cite[Theorem~5.13]{Mukhin:2009et}}]
  \label{thm:mtv-isomorphism}
The map \( \kappa_z\map{\bspec(z)_\mu}{\Shvar(z)_\mu} \) given by \( \kappa_z(\chi) = X_\chi = \ker \chi(\Dd) \) is an isomorphism of schemes. 
\end{Theorem}

For real \( z \) such that \( z_1<z_2< \cdots < z_n \), we obtain a bijection \( \MTV\map{\Shvar(z)_\mu}{\SYT(\mu)} \) by \( \MTV = p\circ\kappa_z^{-1} \).

\section{Speyer's labelling}
\label{sec:speyers-labelling}

A third labelling of \( \Shvar(z)_\mu \) is present in the work of Speyer~\cite{Speyer:2014gg}. The definition is somewhat involved and is described below in as much detail as will be needed later. One important point is that Speyer in fact describes many possible labellings and we choose one that is natural in a particular sense explained later. 

\subsection{Speyer's flat families}
\label{sec:spey-flat-famil}
Let \( \overline{M}_{0,k} \) be the moduli space of stable rational curves with \( k \) marked points. It has a dense open set \( M_{0,k} \) consisting of those curves with a single irreducible component. Fix a curve \( C \in M_{0,k} \) with marked points \( (z_1,z_2,\ldots,z_k) \) and a three element set \( A = \{ 1 \le i_0 < i_1 < i_{\infty} \le k \} \). Let \( \phi_A(C):\PP^1 \longrightarrow \PP^1 \) be the unique isomorphism such that \( \phi_A(C)(z_{i_0})=0 \), \( \phi_A(C)(z_{i_1})=1 \), and \( \phi_A(C)(z_{i_{\infty}}) = \infty \). For each three element set \( A \), we denote a copy of the Grassmanian \( \Gr(r,d)_A =  \Gr(r,d) \). The map \( \phi_A(C) \) induces an isomorphism
\[ \phi_A(C): \Gr(r,d) \longrightarrow \Gr(r,d)_A \]
by \( [x:y] \mapsto \phi_A(C)([x:y]) \). Speyer's family \( \Gg(r,d) \) is the closure of the image of
\begin{equation*}
  \begin{tikzcd}[row sep=0.7]
    M_{0,k} \times \Gr(r,d) \arrow[rr] && \overline{M}_{0,k} \times \prod_{\#A=3} \Gr(r,d)_A \\
    (C,X) \arrow[rr,mapsto] && (C,\phi_A(C)(X)).
  \end{tikzcd}
\end{equation*}
The product runs over all three element subsets of \( [k] \). Given a sequence of partitions \( \blambda = (\lambda^{(1)},\lambda^{(2)},\ldots,\lambda^{(k)}) \) we also define the family \( \Ss(\blambda) = \Gg(r,d)\cap\bigcap_{a \in A}\Shvar(\lambda^{(a)};z_a) \).

\begin{Theorem}[\cite{Speyer:2014gg}]
  \label{thm:speyer-flatness-cm}
The families \( \Gg(r,d) \) and \( \Ss(\blambda) \) are flat and Cohen-Macauley over \( \overline{M}_{0,k} \). Furthermore, if \( \abs{\blambda}=r(d-r) \) then \( \Ss(\blambda)(\RR) \) are a topological covering of \( \overline{M}_{0,k}(\RR) \).
\end{Theorem}

Speyer makes a detailed analysis of the fibres of these families which we summarise here. Given a curve \( C \in \overline{M}_{0,k} \), a \emph{node labelling} for \( C \) is a function \( \nu \) which assigns to every pair \( (C_i,d) \) of an irreducible component of \( C \), and a node \( x \in C_i \), a partition \( \nu(C_i,x) \) with at most \( r \) rows and \( d-r \) columns in such a way that if \( x \in C_i \cap C_j \) then \( \nu(C_i,x)^\comp = \nu(C_j,x) \). Denote the set of node labellings by \( \ttN_C \).

Let \( C_i \subset C \) be an irreducible component. For \( a \in [k] \), let \( d_i(a) \in C_i \) be either the point marked by \( a \) if it is on \( C_i \), or the node by which the marked point is connected to \( C_i \). For a three element set \( A = \{ a_0, a_1, a_\infty \} \subset [k] \) define a map \( \phi_{A,i}\map{\PP^1}{\PP^1} \) that maps \( d_i(a_0) \) to \( 0 \), \( d_i(a_1) \) to \( 1 \) and \( d_i(a_\infty) \) to \( \infty \). We also have an associated map \( \phi_{A,i}\map{\Gr(r,d)}{\Gr(r,d)_A} \). We obtain an embedding \( \Gr(r,d) \hookrightarrow \prod_{\#d_i(A)=3} \Gr(r,d)_A \) given by \( X \mapsto (\phi_{A,i}(X))_{A} \), where the product is over all three element sets \( A \subset [k] \), such that \( \#d_i(A)=3 \). Denote the image by \( \Gr(r,d)_{C_i} \). We will use notation of the form \( \Shvar(\lambda;z)_{C_i} \subset \Gr(r,d)_{C_i}  \).

\begin{Theorem}[{\cite{Speyer:2014gg}}]
  \label{thm:speyer-fibre}
  For a stable curve \( C \in \overline{M}_{0,k} \) with irreducible components \( C_1,C_2,\ldots, C_l \) the fibres of the families \( \Gg(r,d) \) and \( \Ss(\blambda) \) are given by
  \begin{align*}
    \Gg(r,d)(C) &= \bigcup_{\nu \in \ttN_C} \prod_{i} \bigcap_{d \in D_i} \Shvar(\nu(C_i,d);d)_{C_i},  \\
  \Ss(\lambda_{\bullet})(C) &= \bigcup_{\nu \in \ttN_C} \prod_{i} \left( \bigcap_{d \in D_i} \Shvar(\nu(C_i,d),d)_{C_i}  \cap \bigcap_{p \in P_i} \Shvar(\lambda^{(p)},p)_{C_i} \right),
  \end{align*}
  where \( D_i \) is the set of nodes on the component \( C_i \) and \( P_i \) is the set of marked points on \( C_i \).
\end{Theorem}

\subsection{Labelling the fibre}
\label{sec:labelling-fibre}

Now we set \( k=n+1 \) and \( \blambda = (\square^n,\mu^\comp) \) for a partition \( \mu \), of \( n \). Fix \( z \in X_n^< \) and identify this with the stable curve \( C \) with marked points \( z_1,z_2,\ldots,z_n,\infty \) on a single irreducible component. Curves of this type form a connected component \( \mathcal{O} \subset M_{0,n+1}(\RR) \subset \overline{M}_{0,n+1}(\RR) \). We will use the covering \( \Ss(\square^n,\mu)(\RR) \) of \( \overline{M}_{0,n+1}(\RR) \) to label the points of \( \Shvar(z)_\mu \) by \( \SYT(\mu) \). At the point \( C \), the fibre of \( \Ss(\square^n,\mu^\comp) \) is isomorphic to \( \Shvar(z)_\mu \) by construction.

For \( 1 \le q < n \), choose a generic point \( C_q \in \overline{M}_{0,n+1}(\RR) \), in the boundary of this connected component, where the points marked by \( 1,2,\ldots,q \) are on a single irreducible component and the points marked by \( q+1,q+2,\ldots,n,\infty \) on the second component. Now choose any path in \( \overline{\mathcal{O}} \), from \( C \) to \( C_q \) and consider the unique lift of this path to \( \Ss(\square^n,\mu^\comp)(\RR) \) starting at \( X \). By Theorem~\ref{thm:speyer-fibre}, the endpoint of this path over \( C_q \) determines a node labelling \( \nu \) of \( C_q \). Let \( C_q' \) be the irreducible component containing the pint marked by \( \infty \) and \( d \in C_q' \), the unique nodal point. Let \( \mu_q = \nu(C'_q,d) \). Speyer shows that \( \mu_q \subset \mu_{q+1} \) and that \( \abs{\mu_q}=q \). Let \( \Sp(X) \) be the  standard \( \mu \)-tableau determined by the inclusions
\[ \emptyset \subset \mu_1 \subset \mu_2 \subset \cdots \subset \mu_{n-1} \subset \mu. \]
Speyer also shows the resulting map \( \Sp\map{\Shvar(z)_\mu}{\SYT(\mu)} \) is a bijection.

\begin{Remark}
  \label{rem:what-speyer-does}
  In fact Speyer constructs a bijection to a slightly different combinatorial set: the set of \emph{dual equivalence growth diagrams} of a certain shape. These objects are in bijection with standard tableaux, however there is a choice involved that is not entirely natural. The choice we have made above is the only one for which Theorem~\ref{thm:labelling-theorem} is true and it is natural in this sense. 
\end{Remark}

\subsection{Associahedra}
\label{sec:associahedra}

The variety \( \M[k](\RR) \) is a CW-complex and is tiled by associahedra of dimension \( k-3 \). For example, when \( k=5 \), the variety \( \M[5](\RR) \) is tiled by \( 2 \)-associahedra, i.e. by pentagons. For each connected component  \( \Oo \subset \oM(\RR) \), its closure \( \overline{\Oo} \subset \M[k](\RR) \) is such an associahedron. We can lift this CW-structure and the tiling to \( \Ss(\blambda)(\RR) \). The connected components of \( \oM[k](\RR) \), and thus the associahedra tiling \( \M[k](\RR) \) are labelled by circular orderings of \( 1,2,\ldots,k \). If \( \Theta \subset \M[k](\RR) \) is the associahedron corresponding to curves where the points marked by \( 1,2,\ldots,k \) are in increasing order, then we denote by \( \Theta_{pq} \) the facet of \( \Theta \) determined by curves with two irreducible components and where the points marked by \( p,p+1,\ldots,q-1 \) are on one of these components. Since each associahedron is simply connected, the process described in Section~\ref{sec:labelling-fibre} shows that the associahedra in \( \Ss(\square^n,\mu^\comp) \) lying above \( \Theta \subset \M(\RR) \) are labelled by \( \SYT(\mu) \).

\section{Agreement with elementary labelling}
\label{sec:agre-with-elem}

In this section we prove that that the labellings of \( \Shvar(z)_\mu \) described in Sections~\ref{sec:schubert-intersections} and~\ref{sec:speyers-labelling} agree. The proof proceeds by directly checking Speyer's definition is compatible with the limiting process from Section~\ref{sec:an-expl-biject}.

\subsection{Levinson result}
\label{sec:levinson-result}

First, we briefly outline a result of Levinson~\cite{Levinson:2017fp} that will be used in  Section~\ref{sec:coninc-with-spey}.
Recall that \( \Ss(\blambda,\square) \) is a subvariety of \( \overline{M}_{0,n+1}\times\prod_{A} \Gr(r,d)_A \), where \( A \) ranges over three element subsets of \( [n+1] \). Let \( \pi \) be the projection onto \( \prod_{A \subset [n]} \Gr(r,d)_A \). That is \( \pi \) projects onto those Grassmanians for subsets \( A \) which do not contain \( n+1 \).

Suppose \( \abs{\blambda} = r(d-r) - 1 \) and let \( c_{n+1}\map{\M}{\M[n]} \) be the contraction map at the point marked by \( n+1 \) (forgetting the point marked by \( n+1 \) and contracting any unstable irreducible components). The morphism \( c_{n+1} \) allows us to think of \( \Ss(\blambda,\square) \) as a family over \( \M[n](\CC) \) (rather than \( \M(\CC) \)).

\begin{Theorem}[{\cite[Theorem~2.8]{Levinson:2017fp}}]
  \label{thm:Levinson}
  The map \( \pi \) produces an isomorphism onto \( \Ss(\blambda) \) and we have the following commutative diagram,
\begin{center}
  \begin{tikzcd}
    \Ss(\blambda,\square) \arrow[r,"\pi"] \arrow[d] & \Ss(\blambda) \arrow[d] \\
    \M(\CC) \arrow[r,"c_{n+1}"] & \M[n](\CC).
  \end{tikzcd}
\end{center}
Thus \( \pi \) is an isomorphism of families over \( \M[n](\CC) \).
\end{Theorem}

\subsection{The labellings agree}
\label{sec:coninc-with-spey}
We show in this section that the labelling of points described in Section~\ref{sec:an-expl-biject} coincides with Speyer's labelling, that is \( \EL = \Sp \). For clarity of exposition, the bulk of the proof is organised into a series of lemmas below. Let \( X \in \Shvar(z)_\mu \) and let \( T = \Sp(X) \in\SYT(\mu) \). We will use the notation \( X(s), X_\infty \in \Gr(r,d) \) from Section~\ref{sec:an-expl-biject}. We use \( T|_{n-1} \) to denote the tableaux obtained from \( T \) by removing the box containing \( n \). We let \( \Shvar(\square^n,\mu^\comp) \longrightarrow \oM \) be the finite map whose fibre over the point \( z = (z_1,z_2,\ldots,z_n,\infty) \) is \( \Shvar(z)_\mu \). 

\begin{Lemma}
  \label{lem:Xinfty-in-lambda}
Let \( \lambda = \mathrm{sh}(T|_{n-1}) \), then \( X_\infty = \lim_{s \to \infty} X(s) \in \Shvar^\circ(\lambda^\comp;\infty) \).
\end{Lemma}

\begin{proof}
Let \( \iota_\mu \) be the inclusion \( \Shvar(\square^n,\mu^\comp) \hookrightarrow \Ss(\square^n,\mu^\comp) \). That is, for \( (C,X) \)  a point of \( \Shvar(\square^n,\mu^\comp) \subseteq \oM \times \Gr(r,d) \),
\begin{equation*}
  \iota_\mu(C,E) = (\phi_A(C)(E))_{A \subset [n+1]}.
\end{equation*}
If \( p \) is the projection \( \Shvar(\square^n,\mu^\comp) \rightarrow \Gr(r,d) \) and \( p_{\{ 1,2,3 \}} \) is the map \( \Ss(\square^n,\mu^\comp) \rightarrow \Gr(r,d) \) defined by \( p_{\{1,2,3\}}(C,(E_A))=\phi_{\{1,2,3\}}(C)^{-1}(E_{\{1,2,3\}}) \), then we have a commutative diagram
\begin{equation*}
  \begin{tikzcd}
    & \Gr(r,d) \\
    \Shvar(\square^n,\mu^\comp) \arrow[d] \arrow[r,hook,"\iota_\mu"] \arrow[ur,"p"]
    & \Ss(\square^n,\mu^\comp) \arrow[d] \arrow[u,"p_{\{ 1,2,3 \}}"'] \\
    \oM \arrow[r,hook] & \M.
  \end{tikzcd}
\end{equation*}

Let \( C \) be the stable curve with marked points \( z \) and \( C(s) \) the family of stable curves with marked points \( (z_1,\ldots,z_{n-1},s,\infty) \), so \( C = C(z_n) \). We have
\begin{align*}
  X_\infty &= \lim_{s \rightarrow \infty} p(C(s),X(s)) \\
          &= \lim_{s \rightarrow \infty} p_{\{ 1,2,3 \}} \iota_\mu (C(s),X(s)) \\
          &= p_{\{ 1,2,3 \}}Y,
\end{align*}
where \( Y = \lim_{s \rightarrow \infty} \iota_\mu (C(s),X(s)) \).

The point \( Y \) lies over the stable curve \( \lim_{s\to \infty} C(s) \), which has two components, \( C_1 \) with marked points \( z_1,z_2,\ldots,z_{n-1} \) and a node at \( \infty \), and \( C_2 \) with marked points at \( 1 \) and \( \infty \) and a node at \( 0 \). Thus by Theorem~\ref{thm:speyer-fibre},
\begin{equation*}
  Y \in \Shvar(\square^{n-1},\lambda^\comp;z_1,\ldots,z_{n-1},\infty)_{C_1} \times \Shvar(\lambda,\square,\mu^\comp;0,1,\infty)_{C_2}.
\end{equation*}
The partition \( \lambda \) appearing in the node labelling must be \( \mathrm{sh}(T|_{n-1}) \) since we assumed \( T = \Sp(X) \). Now \( p_{\{ 1,2,3 \}} \) is simply projection onto the first factor and is an isomorphism onto \( \Gr(r,d) \) so \( X_\infty = p_{\{ 1,2,3 \}}Y \in \Shvar(\lambda^\comp;\infty) \). However \( X_\infty \notin \Shvar(\nu^\comp;\infty) \) for any \( \nu \) such that \( \abs{\nu} < \abs{\lambda} \) and so we must have that \( X_\infty \in \Shvar^\circ(\lambda^\comp;\infty) \).
\end{proof} 

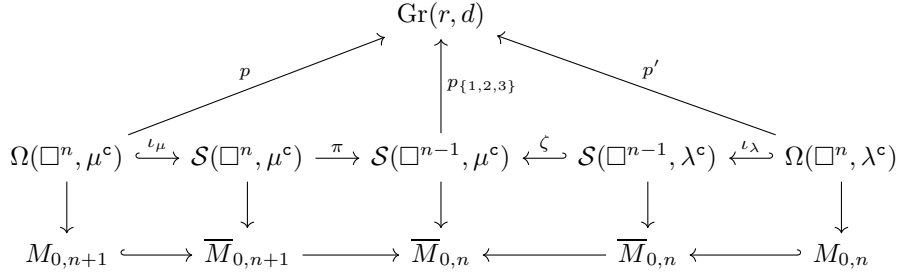
\begin{figure}
  \centering
  \begin{tikzcd}[column sep=17pt]
    & & \Gr(r,d) & & \\ \\
    \Shvar(\square^{n}, \mu^\comp) \arrow[r,hook,"\iota_\mu"] \arrow[d] \arrow[uurr,"p"] & 
        \Ss(\square^{n}, \mu^\comp) \arrow[r,"\pi"] \arrow[d]  &
        \Ss(\square^{n-1},\mu^\comp) \arrow[d] \arrow[uu,"p_{\set{1,2,3}}"'] & 
        \Ss(\square^{n-1},\lambda^\comp) \arrow[l,hook',"\zeta"'] \arrow[d] &
        \Shvar(\square^{n}, \lambda^\comp) \arrow[l,hook',"\iota_\lambda"'] \arrow[d] \arrow[uull,"p'"'] \\
    \oM \arrow[r,hook] & \M \arrow[r] & \M[n] & \M[n] \arrow[l] & \oM[n] \arrow[l,hook',]
  \end{tikzcd}
  \caption{The relationship between various projections}
  \label{fig:big-comm-diagram}
\end{figure}

Our aim will now be to calculate the Speyer labelling of the point \( X_\infty \) in \( \Ss(\square^{n-1},\lambda^\comp) \). However we only have information about the labelling of the points \( X(s) \) in \( \Ss(\square^n,\mu^\comp) \). To relate these two covering spaces we will use Theorem~\ref{thm:Levinson}. With this theorem we produce the large commutative diagram in Figure~\ref{fig:big-comm-diagram}.

In Figure~\ref{fig:big-comm-diagram} \( \iota_\mu \) is the inclusion \( \Shvar(\square^n,\mu^\comp) \hookrightarrow \Ss(\square^n,\mu^\comp) \) described above. The inclusion \( \iota_\lambda \) is defined similarly. The morphism \( \pi \) is the isomorphism appearing in Theorem~\ref{thm:Levinson} and the inclusion \( \zeta \) is induced by the inclusion of \( \Shvar(\lambda^\comp;\infty) \) into \( \Shvar(\mu^\comp;\infty) \).

\begin{Lemma}
  \label{lem:tracing-Xinfty}
Let \( Y = \lim_{s\rightarrow \infty} \iota_\mu(C(s),X(s)) \), and let \( C_\infty \) be the stable curve with marked points \( (z_1,z_2,\ldots,z_{n-1},\infty) \) (the component \( C_1 \) as in the proof of Lemma~\ref{lem:Xinfty-in-lambda}). Then
\begin{equation*}
Y = \pi^{-1}\zeta\iota_{\lambda}(C_\infty,X_\infty).
\end{equation*}
\end{Lemma}

\begin{proof}
  We will show \( p_{\{1,2,3\}}\pi Y = p_{\{1,2,3\}}\zeta\iota_\lambda(C_\infty,X_\infty) \).  Since \( p_{\{1,2,3\}} \) is injective on fibres, this is enough to prove the Lemma.
This  amounts to tracing \( X_\infty \) around the diagram. By commutativity of the diagram,
\begin{align*}
  p_{\{1,2,3\}}\zeta\iota_{\lambda}(C_\infty,X_\infty) &= p'(C_\infty,X_\infty) = X_\infty.
\end{align*}
Now
\begin{align*}
  p_{\{1,2,3\}}\pi Y &= p_{\{1,2,3\}}\pi\lim_{s \rightarrow \infty} \iota_\mu(C(s),X(s)) \\
  &= \lim_{s \rightarrow \infty} p_{\{1,2,3\}}\pi\iota_\mu(C(s),X(s)) \\
  &= \lim_{s \rightarrow \infty} p(C(s),X(s)) = X_\infty. \qedhere
\end{align*}
\end{proof}

\begin{Lemma}
  \label{lem:proj-of-generic-point}
Let \( \Theta \) be the associahedron in \( \Ss(\square^{\abs{\nu}},\nu^\comp) \) labelled by \( S \in \SYT(\nu) \). For generic \( E \in \Theta_{1q} \) 
\begin{equation*}
  p_{\{1,2,3\}}E \in \Shvar^\circ(\tau^\comp;\infty),
\end{equation*}
where \( \tau = \mathrm{sh}(S|_{q}) \).
\end{Lemma}

\begin{proof}
This is a direct application of the  Theorem~\ref{thm:speyer-fibre}, which says if \( E \) is generic then
\begin{equation*}
  E \in \Shvar(\square^q,\tau^\comp;u_1,\ldots,u_q,\infty)_{C_1} \times \Shvar(\tau,\square^{\abs{\nu} - q},\nu^\comp;0,u_{q+1},\ldots,u_{\abs{\nu}},\infty)_{C_2}.
\end{equation*}
Then \( p_{\{1,2,3\}} \) is projection onto the first factor.
\end{proof}

\begin{Remark}
  \label{rem:how-generic-for-prj-lemma}
The proof of Lemma~\ref{lem:proj-of-generic-point} shows in fact we can make the stronger assumption that \( E \) is a generic point of \( \Theta_{1p} \cap \Theta_{1q} \) as long as \( p \geq q \).
\end{Remark}

\begin{Lemma}
  \label{lem:identification-of-action-on-cell}
Let \( \Theta \) be the \( (n-2) \)-associahedron in \( \Ss(\square^n,\mu^\comp) \) labelled by \( T \) and let \( \tilde{\Theta} \) be the \( (n-3) \)-associahedron in \( \Ss(\square^{n-1},\lambda^\comp) \) containing \( \iota_\lambda(X_\infty, C_\infty) \). Then \( \pi^{-1}\zeta(\tilde{\Theta}) = \Theta_{1n} \).
\end{Lemma}

\begin{proof}
Since the maps downstairs in Figure~\ref{fig:big-comm-diagram} are all cell maps, the maps upstairs must also be cell maps. Hence \( \pi^{-1}\zeta(\tilde{\Theta}) \) must be \( \Theta_{ij}' \), the face of some \( (n-2) \)-associahedron \( \Theta' \) in \( \Ss(\square^n,\mu^\comp) \).
Thus \( \Theta_{ij}' \) must contain the point \( \pi^{-1}\zeta\iota_\lambda (C_\infty,X_\infty) \). By Lemma~\ref{lem:tracing-Xinfty}
\begin{equation*}
\pi^{-1}\zeta\iota_\lambda (C_\infty,X_\infty) = Y = \lim_{s \to \infty} \iota_\mu(C(s),X(s)).
\end{equation*}
We know \( \iota_\mu(C(s),X(s)) \in \Theta \) so \( Y \in \Theta_{1n} \). Hence \( \Theta_{ij}' = \Theta_{1n} \).
\end{proof}

\begin{Lemma}
  \label{lem:limit-of-T-point}
\( \Sp(X_{\infty}) = T|_{n-1} \).
\end{Lemma}

\begin{proof}
To show the equality we must show the point \( X_\infty \) is labelled by the tableau \( T|_{n-1} \). That means we must show, for each \( 2 < q < n \) and a generic point \( E \in \tilde{\Theta}_{1q} \), that \( p_{\{1,2,3\}}E \in \Shvar^\circ(\tau^\comp;\infty) \) for \( \tau = \mathrm{sh}(T|_q) \). By commutativity of Figure~\ref{fig:big-comm-diagram}
\begin{equation*}
  p_{\{1,2,3\}}E = p_{\{1,2,3\}}\pi^{-1} \zeta(E).
\end{equation*}
Lemma~\ref{lem:identification-of-action-on-cell} tells us \( \pi^{-1}\zeta(\tilde{\Theta}_{1q}) = \Theta_{1q}\cap\Theta_{1n} \). Since \( X \in \Theta \) and \( \Sp(X)=T \) (which means \( \Theta \) is the \( (n-2) \)-associahedron labelled by \( T \)) and \( \pi^{-1} \zeta(E) \) is generic, we must have that \( p_{\{1,2,3\}}\pi^{-1} \zeta(E) \in \Shvar^\circ(\tau^\comp,\infty) \).
\end{proof}

\begin{Theorem}
  \label{thm:conincides-with-Speyer-label}
We have that \( \Sp = \EL \). That is, if \( X \in \Shvar(\square^n,\mu^\mathtt{c}; z,\infty) \) then the processes described in Sections~\ref{sec:an-expl-biject} and~\ref{sec:labelling-fibre} produce the same tableau.
\end{Theorem}

\begin{proof}
According to Lemma~\ref{lem:Xinfty-in-lambda} \( X_\infty = \lim_{s \rightarrow \infty} X(s) \in \Shvar^\circ(\lambda^\comp;\infty) \). Let \( T = \Sp(X) \). Now Lemma~\ref{lem:limit-of-T-point} tells us \( \Sp(X_\infty) = T|_{n-1} \). This means the tableau \( \EL(X) \) is a tableau of shape \( \mu \) whose restriction to \( n-1 \) is \( T|_{n-1} \). However the unique tableau satisfying these properties is \( T \).
\end{proof}

\section{The algebraic Bethe Ansatz} 
\label{sec:algebr-bethe-ansatz}

Along with Bethe algebras and Schubert intersections, there is a third important player in the story, the critical points of the \emph{master function}. The relationship between these three objects has been studied extensively by Mukhin, Tarasov and Varchenko, see for example~\cite{Mukhin:2012vk}. Critical points have a labelling by standard tableaux in a similar way to points in the spectrum of Bethe algebras (see Section~\ref{sec:spectr-gaud-hamilt}), this is described by Marcus~\cite{Marcus:2010vn} and for the sake of convenience we recall the proof of this result. This will be used to finally identify the MTV and Speyer labellings.

\subsection{Notation}
\label{sec:notation}

Let \( \cartan \) be the Cartan subalgebra of \( \glr \) (so \( \cartan \) is the algebra of diagonal matrices). Let \( \left( \cdot, \cdot \right) \) denote the trace form on \( \glr \) (i.e. the normalised Killing form). Let \( h_i = e_{ii} - e_{i+1,i+1} \) for \( i = 1,\ldots,r-1 \). Let \( \veps_i \in \cartan^* \) be the dual vector to \( e_{ii} \) and \( \alpha_i = \veps_i - \veps_{i+1} \) the dual vector to \( h_i \). With this notation the trace form, transported to \( \cartan^* \), has the following values,
\begin{align}
\label{eq:killing-form}
  \left( \veps_i, \veps_j \right) &= \delta_{ij}, \\
  \left( \alpha_i, \alpha_j \right) &=
                                      \begin{cases}
                                        2  & \text{if } i=j \\
                                        -1 & \text{if } \abs{i-j} = 1 \\
                                        0  & \text{if } \abs{i-j} > 1.
                                      \end{cases}
\end{align}
We identify a partition \( \lambda \) with at most \( r \) parts, with the \( \glr \)-weight \( \sum \lambda^{(i)} \veps_i \).

\subsection{The master function and critical points}
\label{sec:master-function}

Let \( z \in X_n \) be complex parameters, \( \blambda = (\lambda^{(1)},\lambda^{(2)},\ldots,\lambda^{(n)}) \) be a sequence of partitions and let \( \mu \) be a partition such that \( \abs{\mu} = \abs{\lambda_\bullet} \). We also require that there exist non-negative integers \( l_i \) such that \( \mu = \sum_{s} \lambda^{(s)} - \sum_{i=1}^{r-1} l_i \alpha_i \). This last requirement ensures \( \mu \) appears as a weight in \( L(\lamb) \).  We let \( t_i^{(j)} \) be a set of complex variables for \( i = 1,2,\ldots, r-1 \) and \( j = 1,2,\ldots,l_i \).

The \emph{master function} is the rational function
\begin{equation*}
\begin{split}
  & \Phi(\lamb,\mu;z,t) = \Phi(z,t) = \\ &
              \prod_{1 \le a < b \le n} (z_a - z_b)^{\left( \lambda_a,\lambda_b \right)}
              \prod_{a=1}^n \prod_{i = 1}^{r-1} \prod_{j=1}^{l_i} (z_a - t^{(j)}_i)^{-\left( \lambda_a, \alpha_i \right)} 
              \prod_{(i,a) < (j,b)} (t^{(a)}_i - t^{(b)}_j)^{\left( \alpha_i,\alpha_j \right)}.
\end{split}
\end{equation*}
The ordering \( (i,a) < (j,b) \) is taken lexicographically. Let \( S = \log \Phi \). The \emph{Bethe ansatz equations}\index{Bethe ansatz equations} are given by the system of rational functions,
\begin{equation}
\label{eq:bethe-ansatz-equations}
  \frac{\partial S}{\partial t^{(j)}_i} = \frac{\partial}{\partial t^{(j)}_i} \log \Phi(z,t) = 0 \quad \text{for } i = 1,2,\ldots,r-1 \text{ and } j = 1,2,\ldots,l_i.
\end{equation}
A solution to the Bethe ansatz equations is called a \emph{critical point}.\index{critical point} We say a critical point \( t = (t_i^{(j)}) \) is \emph{nondegenerate}\index{critical point!nondegenerate} if the Hessian of \( S \),
\begin{equation*}
  \Hess(S) = \det \left( \frac{\partial S}{\partial t^{(a)}_i \partial t^{(b)}} \right)_{(i,a),(j,b)},
\end{equation*}
evaluated at \( t \) is invertible.

Let \( m = l_1 + l_2 + \ldots + l_{r-1} \). The Bethe ansatz equations are rational functions on \( X_n \times \CC^{m} \), regular away from the finite collection of hyperplanes given by \( t^{(a)}_i - t^{(b)}_j = 0 \). Let \(  \tcrit(\lambda_\bullet)_\mu \) denote the vanishing set of the Bethe ansatz equations, considered as a family over \( X_n \). Let \( \Symcrit \) be the product of symmetric groups \( S_{l_1}\times S_{l_2} \times \cdots S_{l_{r-1}} \subset S_m \), which acts on \( \CC^m \) by permuting the coordinates \( t_i^{(j)} \) with the same lower index. Using~(\ref{eq:killing-form}),
\begin{equation*}
  \prod_{(i,a) < (j,b)} (t^{(a)}_i - t^{(b)}_j)^{\left( \alpha_i,\alpha_j \right)} = \prod_{i=1}^{r-2}\prod_{a=1}^{l_i}\prod_{b=1}^{l_{i+1}} (t_i^{(a)}-t_{i+1}^{(b)})^{-1} \prod_{i=1}^{r-1}\prod_{1\le a<b \le l_i} (t_i^{(a)}-t_i^{(b)})^2.
\end{equation*}
Thus \( \Phi(z,t) \) is invariant under the action of \( \Symcrit \). The quotient \( \tcrit(\lamb)_\mu/\Symcrit \) will be denoted \( \crit(\lamb)_\mu \) and the open subset of nondegenerate critical points \( \crit(\lamb)_\mu^{\nondeg} \). Let \( \infP = \PP\CC[u] \) be the infinite dimensional projective space associated to the polynomial ring. We think of \( \infP \) as the space of monic polynomials. For any \( a \in \NN \), there is an embedding \( \CC^a/S_a \hookrightarrow \infP \) given by sending the orbit of a point \( (t_1,\ldots,t_a) \in \CC^a \) to the unique monic polynomial of degree \( a \), with roots \( t_1,\ldots,t_a \). We will identify \( \crit(\lamb)_\mu \) with its image in \( X_n \times (\infP)^{r-1} \) and denote the tuple of monic polynomials associated to a critical point \( t = (t^{(j)}_i) \) by \( y^t = (y^t_1,\ldots,y^t_{r-1}) \). To be clear, this means if \( t = (t_i^{(j)}) \) is a solution of the Bethe ansatz equations~(\ref{eq:bethe-ansatz-equations}), then \( y^t_i \) is a monic polynomial in \( u \), with roots \( t_i^{(1)}, t_i^{(2)}, \ldots, t_i^{(l_i)} \). Let \( \critfam \) denote the projection \( \crit(\lamb)_\mu \rightarrow X_n \). Denote the fibre of \( \crit(\lamb)_\mu \) over \( z \in X_n \) by \( \crit(\lamb;z)_\mu \).

\begin{Theorem}[{\cite[Theorem~6.1]{Mukhin:2012vk}}]
  \label{thm:weight-function}
There exists a function, called the \emph{universal weight function},\index{universal weight function} \( \omega\map{X_n \times \CC^m/\Symcrit}{L(\lamb)_\mu} \) such that, for a critical point \( t = (t_i^{(j)}) \) in \( \crit(\lamb;z)_\mu \), then
\begin{enumerate}
\item\label{item:omega-singular} \( \omega(z,t) \in L(\lamb)_\mu^\sing \),
\item\label{item:nondeg-nonzero} the critical point \( t \) is nondegenerate if and only if \( \omega(z,t) \) is nonzero, 
\item\label{item:crit-linind} if \( t' \in \crit(\lamb;z)_\mu \) is a critical point distinct from \( t \), and both are nondegenerate then \( \omega(z,t) \) and \( \omega(z,t') \) are linearly independent,
\item\label{item:sim-eigenvector} \( \omega(z,t) \) is a simultaneous eigenvector for \( \uBethe(\lamb;z)_\mu \), and
\item\label{item:eigenvalue-ofGaudin} the eigenvalue of \( H_a(z) \) acting on \( \omega(z,t) \) is
  \begin{equation*}
    \frac{\partial S}{\partial z_a} (z,t). 
  \end{equation*}
\end{enumerate}
\end{Theorem}

\subsection{Examples of the Bethe ansatz equations}
\label{sec:examples-bethe-ansatz}

Below are some examples of the Bethe ansatz equations in simple cases. Explicitly, in full generality, the Bethe ansatz equations are
\begin{equation}
  \label{eq:explicit-BA-eqns}
\begin{split}
\frac{\partial S}{\partial t^{(j)}_i} = 
- \sum_{a = 1}^n \left( \alpha_i,\lambda_a \right) \frac{1}{t^{(j)}_i-z_a} + \sum_{(k,a) \neq (i,j)} \left( \alpha_i,\alpha_k \right) \frac{1}{t^{(j)}_i - t^{(a)}_k} = 0.
\end{split}
\end{equation}

\begin{Example}
  \label{exm:empty-critical-points}
In the case \( n = 1 \), with \( \lamb = (\lambda) \), the only choice for \( \mu \) is \( \mu = \lambda \). Thus \( l_i = 0 \) for all \( i \), the variable \( t \) is simply an empty variable. The master function becomes \( \Phi(z,t) = 1 \). The Bethe ansatz equations in this case are vacuously satisfied and there is a single unique critical point \( t_{\emptyset} \) (the empty critical point). The polynomial \( y^{t_{\emptyset}}_i \) is the unique monic polynomial with no roots, i.e. the constant polynomial \( 1 \). Thus \( \crit(\lambda;z)_\mu \subset (\infP)^{r-1} \) is a single point.
\end{Example}

\begin{Example}
  \label{exm:case-all-box}
In this paper we will be primarily interested in the case \( \lambda_i = \square = \veps_1 \) for all \( i \). In this case \( \abs{\mu} = n \). Since the highest possible weight in \( V^{\otimes n} \) is \( (n) = n\veps_1 \), the integer \( l_i \) is the number of boxes in \( \mu \) sitting strictly below the \( i^{\text{th}} \) row.

In this case \( \left( \veps_1,\veps_1 \right) = 1 \), and \( \left( \alpha_i, \veps_1 \right) = \delta_{1,i} \). Thus the master function becomes
\begin{equation*}
\begin{split}
  \Phi(z,t) = \prod_{1 \le a < b \le n} (z_a - z_b)
              \prod_{a=1}^n \prod_{j=1}^{l_1} (t^{(1)}_j - z_a)^{-1}
              \prod_{i=1}^{r-2}\prod_{a=1}^{l_i}\prod_{b=1}^{l_{i+1}} (t_i^{(a)}-t_{i+1}^{(b)})^{-1} \\
              \prod_{i=1}^{r-1}\prod_{1\le a<b \le l_i} (t_i^{(a)}-t_i^{(b)})^2.
\end{split}
\end{equation*}
\end{Example}

\begin{Example}
  \label{exm:case-n-equals-2}
Next consider the special case when \( n=2 \), so \( \lamb = (\lambda^{(1)},\lambda^{(2)}) \) for some partitions \( \lambda^{(1)} \) and \( \lambda^{(2)} \). Let \( z = (z_1,z_2) \). Make a change of variables
\begin{equation*}
  s^{(j)}_i = \frac{t^{(j)}_i - z_1}{z_2 - z_1}.
\end{equation*}
In these new variables, the Bethe ansatz equations become
\begin{equation*}
\begin{split}
  0  = \frac{\partial S}{\partial s^{(j)}_i} \frac{\partial s^{(j)}_i}{\partial t^{(j)}_i} =  \left( 
  - \frac{(\lambda_1,\alpha_i)}{s^{(j)}_i} - \frac{(\lambda_2,\alpha_i)}{s^{(j)}_i-1} + \sum_{(k,a) \neq (i,j)} \frac{(\alpha_i, \alpha_k)}{s^{(j)}_i - s^{(a)}_k} \right) \frac{1}{z_2-z_1},
\end{split}
\end{equation*}
which can be rearranged to
\begin{equation}
  \label{eq:n=2-bethe-ansatz}
\frac{(\lambda_1,\alpha_i)}{s^{(j)}_i} + \frac{(\lambda_2,\alpha_i)}{s^{(j)}_i-1} = \sum_{(k,a) \neq (i,j)} \frac{(\alpha_i, \alpha_k)}{s^{(j)}_i - s^{(a)}_k},
\end{equation}
and thus do not depend on \( z_1 \) and \( z_2 \). These are the \emph{transformed bethe ansatz equations}.\index{Bethe ansatz equations!transformed} The set of (orbits of) solutions of~(\ref{eq:n=2-bethe-ansatz}) is denoted \( \Scrit(\lambda^{(1)},\lambda^{(2)})_\mu \).

Consider the special case, when \( \lambda_1 = \lambda \) and \( \lambda_2 = \square \). 
By the Pieri rule, for the \( \mu \)-weight space to be nonzero, \( \mu \) must be obtained from \( \lambda \) by adding a single box. Suppose the box is added in row \( e \). Then \( l_i = 1 \) for \( i = 1,2,\ldots,e-1 \) and \( l_i = 0 \) otherwise. Setting \( s_i = s^{(1)}_i \) for \( i=1,2,\ldots, e-1 \), equation~(\ref{eq:n=2-bethe-ansatz}) can be rewritten,
\begin{equation}
  \label{eq:n=2-bethe-ansatz-pieri}
\frac{(\lambda_1,\alpha_i)}{s_i} + \frac{\delta_{1i}}{s_i-1} = \frac{\delta_{i1}-1}{s_i - s_{i-1}} + \frac{\delta_{(i+1)e}-1}{s_i - s_{i+1}}.
\end{equation}
\end{Example}

\begin{Proposition}[{\cite[Lemma~7.2]{Marcus:2010vn}}]
  \label{prp:unique-solution-n2-BAeqn}
There is a unique solution to the transformed Bethe ansatz equations~(\ref{eq:n=2-bethe-ansatz-pieri}), that is \( \Scrit(\lambda,\square)_\mu \) is a single point. In particular 
\begin{equation}
  \label{eq:s1-unique-solution}
  s_1 = 1- \left( \lambda^{(1)} - c \right)^{-1},
\end{equation}
where \( c \) is the content of the box \( \mu \setminus \lambda \).
\end{Proposition}

\subsection{Asymptotics of critical points and Marcus' labelling}
\label{sec:marcus-label-crit}

Later, we will need a result about the asymptotics of critical points as we send the parameters to infinity. Reshetikhin and Varchenko~\cite{Reshetikhin:1995vs} explain how to glue two nondegenerate critical points to obtain a critical point for a larger master function with parameters \( z = (z_1,z_2,\ldots,z_{n+k}) \). Their theorem allows one to track the analytic continuation of this new critical point as we send the parameters \( z_{n+1},\ldots,z_{n+k} \) to infinity and shows that asymptotically we recover the two critical points we started with. The set up for the theorem is the following data, two sequences of partitions,
\begin{itemize}
\item \( \lamb = (\lambda_1,\lambda_2,\ldots,\lambda_n) \), and
\item \( \lamb' = (\lambda_1',\lambda_2',\ldots,\lambda_k') \).
\end{itemize}
Three additional partitions,
\begin{itemize}
\item \( \nu = \sum_i \lambda_i - \sum_j a_j \alpha_j \),
\item \( \nu' = \sum_i \lambda_i' - \sum_j b_j \alpha_j \), and
\item \( \mu = \nu + \nu' - \sum_j c_j \alpha_j = \sum_i \lambda_i + \sum_i \lambda_i' - \sum_j (a_j+b_j+c_j)\alpha_j \),
\end{itemize}
for nonnegative integers \( a_j,b_j \) and \( c_j \). Two nondegenerate critical points
\begin{itemize}
\item \( u = (u^{(j)}_i) \in \crit(\lamb;z)_\nu^\nondeg \), and 
\item \( v = (v^{(j)}_i) \in \crit(\lamb';x)_{\nu'}^\nondeg \),
\end{itemize}
for complex points, \( z = (z_1,z_2,\ldots,z_n) \) and \( x = (x_1,x_2,\ldots,x_k) \); and finally a solution, \( s = (s^{(j)}_i) \in \Scrit(\nu,\nu')_\mu \), to the transformed Bethe ansatz equations~(\ref{eq:n=2-bethe-ansatz}).

\begin{Theorem}[{\cite[Theorem~6.1]{Reshetikhin:1995vs}}]
  \label{thm:RV-colliding-crit-points}
In the limit when \( z_{n+1}, z_{n+2},\ldots, z_{n+k} \) are sent to \( \infty \) in such a way that \( z_{n+i} - z_{n+1} \) remain finite for \( i = 1,2,\ldots,k \), there exists a unique nondegenerate critical point \( t = (t^{(j)}_i) \in \mathcal{C}(\lamb,\lamb';z)_\mu^\nondeg \) such that asymptotically, the critical point has the form
\begin{equation*}
t^{(j)}_i(z) =
\begin{cases}
  u^{(j)}_i(z_1,\ldots,z_n) + O(z_{n+1}^{-1}) &\text{if } 1 \le j \le a_i, \\
  s^{(j)}_i z_{n+1} + O(1)   &\text{if } a_i < j \le a_i + c_i, \\
  v^{(j)}_i(x_1,\ldots,x_k) +z_{n+1} +  O(z_{n+1}^{-1}) &\text{if } a_i + c_i < j \le a_i + b_i + c_i,
\end{cases}
\end{equation*}
where \( x_i = z_{n+i} - z_{n+1} \) for \( i = 1,2,\ldots,k \).
\end{Theorem}

\begin{Corollary}
  \label{cor:RV-theorem-limit-polys}
Let \( t,u,v \) and \( s \) be as in Theorem~\ref{thm:RV-colliding-crit-points}. Taking a limit \( z_{n+i} \to \infty \) such that \( z_{n+i} - z_{n+1} \) is bounded, (which we denote \( \lim_{z\to \infty} \)) we have
\begin{equation*}
  \lim_{z \to \infty} y^t = y^u.
\end{equation*}
\end{Corollary}

\begin{proof}
This is a direct application of Theorem~\ref{thm:RV-colliding-crit-points} to the definition of \( y^t \).
\end{proof}

We restrict our attention to critical points for \( z = (z_1,z_2,\ldots,z_n) \) and \( n \)-tuple of distinct real numbers such that \( z_1 < z_2 < \ldots < z_n \). In the limit when \( z_1, z_2,\ldots,z_n \to \infty \) such that \( z_i = o(z_{i+1}) \), Marcus, \cite{Marcus:2010vn},  describes a method to label critical points in \( \crit(\square^n;z)_\mu \) by standard \( \mu \)-tableaux. Marcus' theorem is recalled below, along with the proof. Recall, if \( T \in \SYT(\mu) \) then \( T|_{n-1} \) is the tableaux obtained by removing the box containing \( n \) from \( T \).

\begin{Theorem}[{\cite[Theorem~6.1]{Marcus:2010vn}}]
  \label{thm:Marcus-main-theorem}\index{tableau!standard}
Given a standard tableaux, \( T \), of shape \( \mu \), there is a unique critical point \( t_T \in \crit(\square^n;z)_\mu \) such that, if we set \( y^T = y^{t_T} \), then
\begin{equation}
\label{eq:lim-zn-restrict-tab}
  \lim_{z_n \to \infty} y^T = y^{T|_{n-1}},
\end{equation}
and taking the limit \( z_1,z_2,\ldots,z_n \) such that \( \abs{z_i} << \abs{z_{i+1}} \), asymptotically
\begin{equation}
\label{eq:limit-eigenvalues}
  z_a \frac{\partial S}{\partial z_a} \sim c_T(i) + O(z_i^{-1}).
\end{equation}
\end{Theorem}

\begin{proof}
We will prove this by induction on \( n \). For \( n=1 \), the only partition is \( \square \) and thus there is a unique tableau \( T \). From Example~\ref{exm:empty-critical-points} we know \( \crit(\square;z)_\square \) contains a unique critical point, the empty critical point \( t_\emptyset \) and we simply set \( t_T = t_\emptyset \). Thus \( y^T = 1 \). The equations~(\ref{eq:lim-zn-restrict-tab}) and~(\ref{eq:limit-eigenvalues}) are vacuously satisfied.

For general \( n \), we will use Theorem~\ref{thm:RV-colliding-crit-points} to inductively build a critical point corresponding to \( T \in \SYT(\mu) \). Let \( \lambda = \mathrm{sh}(T|_{n-1}) \), the partition obtained by removing the box labelled \( n \) in \( T \), from \( \mu \). By induction, there is a unique critical point \( t_{T|_{n-1}} \in \crit(\square^{n-1};z_1,\ldots,z_{n-1})_\lambda \). To build a critical point in \( \crit(\square^n;z)_\mu \), we need to fix a critical point in \( \crit(\square;0)_\square \), and a transformed critical point in \( \Scrit(\lambda,\square)_\mu \). The former contains only the empty critical point and the latter contains a unique point \( s = (s_1,\ldots,s_{e-1}) \) (where \( e \) is the row containing \( n \) in \( T \)) by Proposition~\ref{prp:unique-solution-n2-BAeqn}, where
\begin{equation}
  \label{eq:s1-content-n}
s_1 = 1- \left( \lambda^{(1)} - c_T(n) \right)^{-1}.
\end{equation}
Thus, given \( t_{T|_{n-1}} \), and the data of where to add an \( n^{\text{th}} \) box to \( T|_{n-1} \), we obtain by Theorem~\ref{thm:RV-colliding-crit-points} a unique critical point \( t_T \in \crit(\square^n;z)_\mu \). By Corollary~\ref{cor:RV-theorem-limit-polys} we obtain~(\ref{eq:lim-zn-restrict-tab}).

All that is left to do is to prove~(\ref{eq:limit-eigenvalues}). This will also be done by induction. We need to investigate the eigenvalues
\begin{equation*}
  z_a \frac{\partial S}{\partial z_a} = \sum_{b \neq a} \frac{z_a}{z_a - z_b} - \sum_{j=1}^{\abs{\mu} - \mu^{(1)}} \frac{z_a}{z_a - t^{(j)}_1},
\end{equation*}
of the operators \( z_a H_a(z) \). Suppose first that \( a < n \), then
\begin{equation*}
  z_a \frac{\partial S}{\partial z_a} = \sum_{\substack{b \neq a \\ b < n}} \frac{z_a}{z_a - z_b} - \sum_{j=1}^{\abs{\mu} - \mu^{(1)}-\delta_{e>1}} \frac{z_a}{z_a - t^{(j)}_1} + \frac{z_a}{z_a-z_n} - \frac{\delta_{e>1}z_a}{z_a - t^{(\abs{\mu}-\mu^{(1)})}_1}.
\end{equation*}
But in the limit \( z_a/(z_a-z_n) \sim 0 \) and by Theorem~\ref{thm:RV-colliding-crit-points} \( t^{(\abs{\mu}-\mu^{(1)})}_1 \sim s_1 z_{n} + O(1)  \) so
\begin{align*}
  z_a \frac{\partial S}{\partial z_a} &\sim \sum_{\substack{b \neq a \\ b < n}} \frac{z_a}{z_a - z_b} - \sum_{j=1}^{\abs{\mu} - \mu^{(1)}-\delta_{e>1}} \frac{z_a}{z_a - t^{(j)}_1} \\
&= \sum_{\substack{b \neq a \\ b < n}} \frac{z_a}{z_a - z_b} - \sum_{j=1}^{\abs{\lambda} - \lambda^{(1)}} \frac{z_a}{z_a - t^{(j)}_1} \\
&= z_a \frac{\partial S'}{\partial z_a},
\end{align*}
where \( S' = S(\square^{n-1},\lambda;z_1,\ldots,z_{n-1}) \) is the logarithm of the master function for the weight \( \lambda \). Thus by induction~(\ref{eq:limit-eigenvalues}) holds for \( a < n \).

Now we need to check~(\ref{eq:limit-eigenvalues}) for \( a=n \). This turns out to be a simple calculation. By Theorem~\ref{thm:RV-colliding-crit-points}
\begin{align*}
  z_n \frac{\partial S}{\partial z_n} &= \sum_{b < n} \frac{z_n}{z_n - z_b} - \sum_{j=1}^{\abs{\mu} - \mu^{(1)}} \frac{z_n}{z_n - t^{(j)}_1} \\
&\sim (n-1) - (\abs{\mu} - \mu^{(1)} - \delta_{e>1}) - \frac{\delta_{e>1}}{1-s_1} + O(z_n^{-1}) \\
&= \mu^{(1)} - \delta_{e1} - \frac{\delta_{e>1}}{1-s_1} + O(z_n^{-1}).
\intertext{Using~(\ref{eq:s1-content-n}),}
  z_n \frac{\partial S}{\partial z_n} &\sim \mu^{(1)} - \delta_{e1} - \delta_{e>1} \left( \lambda^{(1)} - c_T(n) \right) + O(z_n^{-1}).
\end{align*}
If \( e > 1 \) then \( \mu^{(1)} = \lambda^{(1)} \), and if \( e = 1 \) then \( c_T(n) = \mu^{(1)} - 1 \) so the Theorem is proved.
\end{proof}

\section{Critical points and Schubert intersections}
\label{sec:crit-points-schub}

In this final section we describe a relationship between critical points and Schubert intersections called the coordinate map. We use this and the fact that critical points determine the eigenvalues of the Bethe algebras to prove that the MTV and Speyer labellings agree. The most important part of the argument is a careful analysis of exactly when the coordinate map is continuous.

\subsection{The coordinate map}
\label{sec:second-mtv-iso}

This section describes the relationship between critical points for the master function and Schubert intersections. 
Let \( X \in \Gr(r,d) \). Since \( \Gr(r,d) \) is paved by Schubert cells \( X \in \Shcell(\mu^\comp;\infty) \) for some unique partition \( \mu \). By definition \( X \) is an \( r \)-dimensional vector space of polynomials in the variable \( u \), of degree less than \( d \). Let \( l_i \) be the number of boxes below the \( i^{\text{th}} \) row in \( \mu \) (c.f Example~\ref{exm:case-all-box}). Set \( d_i = \mu_i + r -1 \). We can choose a ordered basis \( f_1(u), f_2(u), \ldots, f_r(u) \) of monic polynomials with descending degrees \( d_i \). Consider the polynomials 
\begin{equation*}
  y_a(u) = \mathrm{Wr}(f_{a+1}(u),\ldots,f_r(u)), \quad a = 0,1,\ldots,r-1.
\end{equation*}
The polynomial \( y_a \) has degree \( l_a \). Denote its roots by \( t^{(1)}_a,t^{(2)}_a,\ldots,t^{(l_a)}_{a} \). The polynomial \( y_a(u) \) determines a point in \( \infP \). The following lemma demonstrates the polynomials \( y_a(u) \in \infP \) depend only on \( X \) and not the basis chosen.

\begin{Lemma}
  \label{lem:MTV-morphism-well-defined}
Suppose \( \{ f_i(u) \} \) and \( \{ f'_i(u) \} \) are two bases of \( X \) of monic polynomials with descending degrees. Then 
\begin{equation*}
  \Wr(f_{a+1}(u),\ldots,f_r(u)) = \alpha \Wr(f'_{a+1}(u),\ldots,f'_r(u))
\end{equation*}
for some scalar \( \alpha \in \CC \).
\end{Lemma}

\begin{proof}
We use the fact that the descending sequence \( d_1 > d_2 > \ldots >d_r \)  of degrees for any such basis is determined entirely by the partition \( \mu \). That is, \( \deg f_i(u) = \deg f'_i(u) = d_i \). By Lemma~\ref{lem:linear-class-wronskian} the Wronskian \( \Wr(g_1,\ldots,g_k) \) is determined by the space spanned by the polynomials \(  g_1,\ldots,g_k \), we must prove
\begin{equation}
\label{eq:span-f-equals-span-fprime}
  \CC\set{ f_{a+1}(u), \ldots, f_r(u) } = \CC\set{ f'_{a+1}(u), \ldots, f'_r(u) }.
\end{equation}
Since both bases span \( X \), \( f'_a(u) = \alpha_1 f_1(u) + \alpha_2 f_2(u) + \ldots + \alpha_r f_r(u) \) for some complex numbers \( \alpha_i \). But the degrees of the \( f_i \) are strictly descending so \( \alpha_i = 0 \) for \( i > a \). Hence \( f'_a(u) \in \CC\{ f_a(u), \ldots, f_r(u) \} \). By induction (\ref{eq:span-f-equals-span-fprime}) must be true.
\end{proof}

The map \( \mtvSC\map{\Gr(r,d)}{\left( \infP \right)^r} \) defined by 
\begin{equation*}
  \mtvSC(X) = (y_a)_{a=0}^{r-1} =  \left( \Wr(f_{a+1}(u),\ldots,f_r(u)) \right)_{a=0}^{r-1}
\end{equation*}
for some choice of monic basis of descending degrees \( f_1(u),f_2(u),\ldots,f_r(u) \), of \( X \) is called the \emph{coordinate map}.

\begin{Remark}
  \label{rem:MTV-map-not-cont}\index{coordinate map!continuity of}
The coordinate map is not continuous! This is easily seen in an example. Let \( r=2 \) and \( d=3 \). Consider the \( 1 \)-parameter family of subspaces
\begin{equation*}
  X(s) = \CC \{ u^2 + s, u \}.
\end{equation*}
In this case \( \mtvSC(X(s)) = \left( s-u^2, u \right) \). However
\begin{equation*}
   X_\infty = \lim_{s \to \infty}X(s) = \CC\{ 1, u \}.
\end{equation*}
So \( \mtvSC(X_\infty) = \left( 1, 1 \right) \), which is clearly not the same as \( \lim_{s \to \infty} \mtvSC(X(s)) = \left( 1, u \right) \).
\end{Remark}

Essentially the problem in Remark~\ref{rem:MTV-map-not-cont} is the monic basis of descending degrees which we are using to calculate \( \mtvSC(X(s)) \) no longer has descending degrees in the limit. Whenever we can find a continuous family of monic bases of descending degrees the map \( \mtvSC \) will be continuous  since taking the Wronskian of a tuple of polynomials is algebraic. An important case in which we can do this is for a a Schubert cell. If \( X \in \Shcell(\mu^\comp;\infty) \), we can find a unique basis of the form 
\begin{equation*}
  f_i(u) = u^{d_i} + \sum_{\substack{j = 1 \\ d_i-j \notin \mathbf{d}}} a_{ij} u^{d_i - j},
\end{equation*}
where \( \mathbf{d} = (d_1,d_2,\ldots,d_r) \). The \( a_{ij} \) define algebraic coordinates on \( \Shcell(\mu^\comp;\infty) \). Hence \( \mtvSC \) is algebraic when restricted to any open Schubert cell. In Section~\ref{sec:partial-continuity} we will prove that the coordinate map is continuous along certain paths in \( \Gr(r,d) \) which are allowed to have limit points outside a Schubert cell.


\begin{Theorem}[{\cite[Theorem~5.3]{Mukhin:2012vk}}]
  \label{thm:coord-map-onto-crit}
The image of \( \Shcell(\square^n,\mu^\comp;z,\infty) \) under the coordinate map is contained in \( \crit(\square^n;z)_\mu \).
\end{Theorem}

\subsection{Partial continuity of the coordinate map}
\label{sec:partial-continuity}

\index{coordinate map}
\index{coordinate map!continuity of}
We will need a little more information about when the coordinate map is continuous. Let \( \mu \) be a partition and let \( \lambda \subseteq \mu \) be a partition with one less box, that is \( \abs{\lambda} = \abs{\mu} -1 \). Denote by \( e \) the row of \( \mu \) from which we need to remove a box to obtain \( \lambda \). 

Let \( d_i = \mu_i + r - 1  \) and \( d'_i = \lambda_i +r -1 \). We denote the respective decreasing sequences by \( \mathbf{d} = (d_1,d_2,\ldots,d_r) \) and \( \mathbf{d}' = (d'_1,d'_2,\ldots,d'_r) \). Recall that \( X \in \Shcell(\mu^\comp;\infty) \) (respectively \( X \in \Shcell(\lambda^\comp;\infty) \)) if and only if there exists a basis \( f_1,f_2,\ldots,f_r \) of \( X \) such that \( \deg f_i = d_i \) (respectively \( \deg f_i = d'_i \)). In particular if \( f \in X \) then \( \deg f \in \mathbf{d} \) (respectively \( \deg f \in \mathbf{d}' \)). Since we removed a single box from \( \mu \) in row \( e \) to obtain \( \lambda \) we have
\begin{equation*}
  d'_i =
  \begin{cases}
    d_i & \text{if } i \neq e \\
    d_e - 1 & \text{if } i = e.
  \end{cases}
\end{equation*}

Fix \( X \in \Shcell(\mu^\comp;\infty) \). Let \( X(s) \in \Shcell(\mu^\comp;\infty) \) be a continuous one-parameter family over \( \RR \) of subspaces. We thus have a unique basis
\begin{equation*}
  f_i(u;s) = u^{d_i} + \sum_{\substack{j = 1 \\ d_i-j \notin \mathbf{d}}} a_{ij}(s) u^{d_i - j},
\end{equation*}
for \( X(s) \), for each \( s \in \RR \). The \( a_{ij}\map{\RR}{\CC} \) are continuous functions. 

\begin{Lemma}
  \label{lem:bounds-on-aij}
The limit point \( X_\infty \) of this family is contained in  \( \Shcell(\lambda^\comp,\infty) \) if and only if
\begin{align}
  \label{eq:lim_a_ij-i-not=-e}
  \lim_{s \to \infty} \abs{a_{ij}(s)} &< \infty \quad \text{for } i \neq e \\
\intertext{and}
  \label{eq:lim_a_e1}
  \lim_{s \to \infty} \abs{a_{e1}(s)} &= \infty, \\
  \label{eq:lim_a_ej-div-a_e1}
  \lim_{s \to \infty} \abs{\frac{a_{ej}(s)}{a_{e1}(s)}} &< \infty.
\end{align}
\end{Lemma}

\begin{proof}
If properties (\ref{eq:lim_a_ij-i-not=-e}), (\ref{eq:lim_a_e1}) and (\ref{eq:lim_a_ej-div-a_e1}) hold then in the limit, the basis \( f_1,\ldots,f_r \) is a sequence of \( r \) polynomials which have descending degrees \( d'_1 > d'_2 > \ldots > d'_r \). Hence \( X_\infty \in \Shcell(\lambda^\comp;\infty) \).

In the other direction, if (\ref{eq:lim_a_ij-i-not=-e}) fails to hold, then there exists an \( i \neq e \) and a \( j \) such that \( \lim_{s \to \infty} \abs{a_{ij}(s)} = \infty \). There are two cases. First if \( i < e \), then choose \( j \) such that \( a_{ij}(s) \) has the fastest growth (for fixed \( i \)). Thus the polynomial
\begin{equation*}
  \lim_{s \to \infty} a_{ij}(s)^{-1}f_i(u;s) \in X_\infty
\end{equation*}
has degree \( d_i-j \). But \( d_i-j \notin \mathbf{d} \) and \( d_i-j < d_e-1 \) thus \( d_i-j \notin \mathbf{d}' \) and so we must have \( X_\infty \notin \Shcell(\lambda^\comp;\infty) \).

Now for the second case, assume there exists an \( i \) such that \( i > e \) and such that \( \lim_{s \to \infty} \abs{a_{ij}(s)} = \infty \). For any \( f \in X_\infty \) there exist functions \( \alpha_i\map{\RR}{\CC} \) such that
\begin{equation*}
  f = \lim_{u \to \infty} \alpha_1(s) f_1(u;s) + \alpha_2(s) f_2(u;s) + \ldots + \alpha_r(s) f_r(u;s).
\end{equation*}
If \( X_\infty \in \Shcell(\lambda^\comp;\infty) \) then we can choose \( f \) so that \( \deg f = d_i \). In order for this to be true we first need \( \lim_{u \to \infty} \alpha_k(s) = 0 \) for \( k > i \) (since \( \lim_{u \to \infty} f_k(u;s) \) exists in this case and has degree \( d_k > d_i \)). By assumption
\begin{equation*}
  \deg \lim_{u \to \infty} \alpha_i(s) f_i(u;s) < d_i
\end{equation*}
and for \( k < i \) we also have 
\begin{equation*}
  \deg \lim_{u \to \infty} \alpha_k(s) f_k(u;s) < d_i.
\end{equation*}
Thus we have a contradiction and \( X_\infty \notin \Shcell(\lambda^\comp;\infty) \).

Finally if conditions (\ref{eq:lim_a_e1}) and (\ref{eq:lim_a_ej-div-a_e1}) do not hold then let \( j \) be minimal with respect to the condition \( a_{ej}(s) \) has the largest order of growth. By assumption \( j > 1 \). Then the degree of \( \lim_{u \to \infty} a_{ej}(s)^{-1}f_e(u;s) \) is \( d_e-j \notin \mathbf{d}' \). Hence \( X_\infty \notin \Shcell(\lambda^\comp;\infty) \).
\end{proof}

Using this lemma we can prove the continuity of the coordinate map \( \mtvSC \) along certain kinds of paths. Heuristically, the paths along which \( \mtvSC \) is continuous are those in \( \Gr(r,d) \) which remain inside a Schubert cell and if they pass into another Schubert cell, do so only in a way in which the partition labelling the Schubert cell is a single box smaller than the partition labelling the original Schubert cell. Let \( X, X(s) \) and \( X_\infty \) be as above.

\begin{Proposition}
\label{prp:cont-of-theta-along-paths}
If \( X_\infty \in \Shcell(\lambda^\comp;\infty) \) then 
\begin{equation*}
  \mtvSC(X_\infty) = \mtvSC\left( \lim_{u \to \infty} X(s) \right) = \lim_{u \to \infty} \mtvSC(X(s)).
\end{equation*}
\end{Proposition}

\begin{proof}
Since \( X_\infty \in \Shcell(\lambda^\comp;\infty) \) by Lemma~\ref{lem:bounds-on-aij} we have conditions  (\ref{eq:lim_a_ij-i-not=-e}), (\ref{eq:lim_a_e1}) and (\ref{eq:lim_a_ej-div-a_e1}) and thus have a monic basis of descending degrees
\begin{align*}
  f_i^\infty(u) = \lim_{u \to \infty} f_i(u;s) &= u^{d_i} + \sum_{\substack{j = 1 \\ d_i-j \notin \mathbf{d}}} a_{ij}^\infty u^{d_i - j}, \quad \text{for } i \neq e, \\
  f_e^\infty(u) = \lim_{u \to \infty} a_{e1}(s)^{-1} f_e(u;s) &= u^{d_e-1} + \sum_{\substack{j = 2 \\ d_i-j \notin \mathbf{d}}} b_{j}^\infty u^{d_i - j}.
\end{align*}
Here \( a^\infty_{ij} = \lim_{u \to \infty}a_{ij}(s) \) and \( b_{j}^\infty = \lim_{u\to\infty} a_{ej}(s)/a_{e1}(s) \) which exist by  (\ref{eq:lim_a_ij-i-not=-e}) and (\ref{eq:lim_a_ej-div-a_e1}).

Let \( X^a(s) = \CC\{ f_a,\ldots,f_r \} \) and \( X^a_\infty(s) = \CC\{ f^\infty_a,\ldots,f^\infty_r \} \). We can use these spaces to calculate the Wronskian. That is
\begin{align*}
  \mtvSC(X(s)) &= \left( \Wr(X^{1}(s)), \ldots, \Wr(X^r(s)) \right) \\
\intertext{and}
  \mtvSC(X_\infty) &= \left( \Wr(X^{1}_\infty(s)), \ldots, \Wr(X^r_\infty(s)) \right)
\end{align*}
Since \( \lim_{s \to \infty}X^a(s) = X^a_\infty(s) \) and since the Wronskian is continuous 
\begin{align*}
  \lim_{u \to \infty} \mtvSC(X(s)) &= \lim_{s \to \infty} \left( \Wr(X^{1}(s)), \ldots, \Wr(X^r(s)) \right) \\
  &= \left( \Wr(\lim_{s \to \infty} X^{1}(s)), \ldots, \Wr(\lim_{s \to \infty} X^r(s)) \right) \\
  &=  \left( \Wr(X^{1}_\infty(s)), \ldots, \Wr(X^r_\infty(s)) \right) \\
  &= \mtvSC(X_\infty). \qedhere
\end{align*}
\end{proof}

\subsection{Identifying the Speyer and MTV labellings}
\label{sec:ident-spey-mtv}

We are now in a position to prove the main theorem, that the Speyer and MTV labellings agree. In Theorem~\ref{thm:conincides-with-Speyer-label} we saw that \( \Sp = \EL \) so it will be enough to show that \( \MTV = \EL \). First we show that Speyer's labelling is compatible with Marcus' and the coordinate map.

\begin{Theorem}
  \label{thm:speyer-equals-marcus}
If \( X \in \Shvar(z)_\mu \) and \( \EL(X)=T \in \SYT(\mu) \), then \( \theta(X) = t_T \). 
\end{Theorem}

\begin{proof}
  We prove this theorem by induction on \( n \). In the case \( n = 1 \), Both \( \Shvar(z)_\square \) and \( \crit(\square;z)_\square \) consist of a single point both of which are labelled by the unique \( \square \)-tableaux. Thus by Theorem~\ref{thm:coord-map-onto-crit} they are mapped to each other by \( \theta \).

For \( n > 1 \), as above, let \( X(s) \in \Shvar(\square^n,\mu^\comp;z_1,\ldots,z_{n-1},s,\infty) \) be the unique family of points passing though \( X \). 
By Lemma~\ref{lem:Xinfty-in-lambda}, the limit \( X_\infty = \lim_{s \to \infty}X(s)  \) is contained in \( \Shcell(\lambda^\comp;\infty) \) where \( \lambda = \mathrm{sh}(T|_{n-1}) \). We also know by Lemma~\ref{lem:limit-of-T-point} that \( \EL(X_\infty) = T|_{n-1} \). By the induction hypothesis \( \mtvSC(X_\infty) = t_{T|_{n-1}} \). In particular, by Proposition~\ref{prp:cont-of-theta-along-paths}
\begin{equation*}
  \lim_{s \to \infty} \mtvSC(X(s)) = \mtvSC(X_\infty) = t_{T|_{n-1}}.
\end{equation*}
However, by Theorem~\ref{thm:RV-colliding-crit-points} and by the definition of Marcus' labelling there is a unique family of critical points with this property, namely the family passing through \( t_T \). Hence we have that \( \mtvSC(X_T) = t_T \).
\end{proof}

\begin{Theorem}
  \label{thm:MTV-equals-Speyer}
  If \( X \in \Shvar(z)_\mu \) then \( \MTV(X) = \EL(X) = \Sp(X) \).
\end{Theorem}

\begin{proof}
  We recall briefly from Section~\ref{sec:mtv-labelling} how \( \MTV\map{\Shvar(z)_\mu}{\SYT(\mu)} \) is defined. We consider the functional \( \chi = \kappa_z^{-1}(X) \in \bspec(z)_\mu \) and then take a limit as \( z_i \to \infty \) such that \( z_i/z_{i+1} \to 0 \). The eigenvalues \( \lim_{z\to\infty} \chi(z_aH_a) \) determine the content of the box containing \( a \) in the tableau \( \MTV(X) \).

  By \cite[Corollary~8.7]{Mukhin:2012vk} \( \omega \circ \mtvSC(X) \) is a simultaneous eigenvector for \( \uBethe(z)_\mu \) with eigenvalues given by \( \chi = \kappa_z^{-1}(X) \in \bspec(z)_\mu \). Let \( \EL(X)=T \in \SYT(\mu) \). Then by Theorem~\ref{thm:speyer-equals-marcus} \( \theta(X) = t_{T} \) and by Theorem~\ref{thm:weight-function},~(\ref{item:eigenvalue-ofGaudin}),
\begin{equation*}
  \chi(z_a H_a(z)) = z_a \frac{\partial S}{\partial z_a}(z,t_T),
\end{equation*}
so by Theorem~\ref{thm:Marcus-main-theorem},
\begin{equation*}
  \chi(z_a H_a(z)) = c_T(a) + O(z_a^{-1}).
\end{equation*}
which implies that \( \MTV(X) = T = \EL(T) \).
\end{proof}

\bibliographystyle{alpha}
\bibliography{bibliography}

\end{document}